# SPIRAL CHAINS: A NEW PROOF OF THE FOUR COLOR THEOREM

**I. Cahit**[1]

**Abstract**

Acceptable but due to extensive usage of a computer rather unpleasant proof of the famous four color map problem of Francis Guthrie were settled eventually by W. Appel and K. Haken in 1976. Using the same method but shortening the proof twenty years later by another team, namely N. Robertson, D. P. Sanders, P. D. Seymour and R. Thomas would not improve considerably the readability of the proof either. Thus it has been widely accepted the need of more elegant and readable proof. There are considerable number of equivalent formulations of the problem but none of them promising for a possible non-computer proofs. On the other hand known proofs are used the concept of Kempe chain and reducibility of the configurations which were a century old ideas. With these in mind we have introduced a new concept which we call "spiral chains" in the maximal planar graphs and instead of sticking to the reducible configurations which were the main complexity issue in the classical proof, we prefer to show the four colorability of the spiral chains by starting coloring from the inner vertices of a maximal planar graphs. We have shown that for any maximal planar graph as long as spiral chains are being used we do not need the fifth color. Henceforth this paper offers another proof to the four color theorem which is not based on deep and abstract theories from the other branches of mathematics or using computing power of computers, but rather completely on a new idea in graph theory.

## 1 Introduction

The famous four color problem is to color the regions of a given map by using least number of colors so that adjacent regions have different colors. It is easy to find maps that require more than three colors and quite early result five colors is enough for all maps. But when the number of color is four, it took almost a century for a correct proof. Unfortunately the proof technique is heavily based on the use of a computer and in fact one may speculate that if it were electronic computers in the 18$^{th}$ century when we would not probably be waiting too long for a proof. Of course assuming the same proof technique basically is being used. Actually the first refuted proof of Kempe was the main inspiration in the Appel and Haken proof, that is assume that all but one region of the given map properly colored by four colors. Then the notion of reducibility, charging-discharging are involved in the proof. In our proof we will stick to that the four color map problem is mainly a problem in topological graph theory since there are many attempts to pull the problem on different mathematical grounds. In fact in the new proof we show that all maximal planar graphs are four colorable as in the previous proofs but without using the notion of reducibility. Of course the notion of reducibility is useful not only in the proof of the four color theorem but in the others [16] but there is a reason as not to stick on only at the reducibility as Ekhad and Zeilberger quotes [15]:

> *A celebrated example of an a posteriori trivial theorem is Appel and Haken's Four-Color Theorem [2]. Their approach was to find a finite unavoidable set of reducible configurations. The original proof used an excessive amount of human efforts. This was considerably improved in the new proof by Robertson, Sanders, Seymour and Thomas*

[1] icahit@ebim.com.tr

*[3],and this is not the ultimate proof. Eventually one should be able to type* `Prove 4CT( );` *and the truth of the theorem should be implied by the halting of the program. In order check the validity, the checker would not need to see any specific configuration. Everything should be done internally and silently by the computer. All that the checker would have to do is check the program.*

With these in mind we have introduced a new concept which we called "spiral chains" in the maximal planar graphs and instead of sticking to the reducible configurations which were the main complexity issue in the classical proofs, we prefer to show the four colorability of the spiral chains starting the inner vertices of a maximal planar graphs. We have shown that for any maximal planar graph as long as spiral chains are being used we do not need the fifth color. Henceforth this paper offers another proof to the four color theorem which is not based on a deep and abstract theories from the other branch of mathematics or power of computers, but rather completely on a new idea in graph theory which we call *spiral chains* of the maximal planar graphs.
In a way our proposition of using the spiral chains of a maximal planar graphs in the proof of the four-color theorem is what is meant by the *ultimate* proof mentioned above. Clearly we will have to prove that our algorithm will halt.

# 2 Can Kempe-chain method be repaired?

Kempe chain was a property of proper coloring as well as a tool to re-establish proper coloring e.g., by exchanging consecutive colors over the chain. As it has been shown by Heawood and many others that this the last minute life-save operation would not always possible. In recent works Wagon [7] and Hutchinson and Wagon [9] have investigated Kempe failure of Kempe chain in the famous *bad*[2] examples and actually Wagon raised the claim that with a very high probability re-coloring the whole graph (or giving a second chance) there would be no *impasse* for the Kempe method. But even if Wagon's conjecture were would be proved to be true its conclusion would be something that "all planar graphs are four colorable with probability close to 1". More specifically Wagon asked the following:

- *Question:* Is it true that any planar graph admits a labeling for which Kempe's coloring method works? That is, does the algorithm- the "relabel and retry" strategy- always halt?

Then he raised a conjecture with J. Hutchinson that "*The Kempe-Kittell algorithm always halts*". For completeness we give the algorithm below.

***Algorithm:*** The Kempe-Kittell Algorithm.
1. Given a planar graph, label the vertices randomly.
2. Set the Kempe order of the vertices, by removing a vertex of degree 5 or less, and repeating until only a single vertex is left.
3. Color the single vertex red.
4. Add vertices back in the reverse order in which they were removed, coloring them as follow:
   I. Follow Kempe's method to letter, switching color on chains defined by nonconsecutive vertices in an attempt to free up a color.

[2] Counterexample in one sense.

II. If Kempe's method fails use Kittell's random Kempe-chain switches until the impasse is resolved.

The above idea is nice but it would be very difficult to prove it for any planar graph since it contains two *random* processes i.e., use of random vertex coloring and use of random Kittell's Kempe-chain switching. That is from the discussions above we conclude that we know how and why sometimes Kempe method fails but unfortunately we cannot be sure when this would be happen. And because of this it would be extremely difficult to give the halting conditions for the Kempe-Kittell algorithm.

Another approach for map coloring, under the context of AI has been given by A. Parmar [14]. She has investigated map coloring problem under mathematical structures antimatroids, also known as shelling structures, a construct used to formalize when greedy (local) algorithms are optimal, as well as their relation to the strong measure of progress P. Map coloring is one domain where planning can be reduced to reasoning about the order in which goals are achieved. The Four Color Theorem (Appel & Haken 1977a; 1977b [2]) guarantees that any map can be colored with at most four colors, and there is already a quadratic algorithm (Robertson *et al.* 1995 [3]) for doing so. Her rationale for researching this domain is knowing when and how we can solve the problem *without having to backtrack*. This is the same motivation as that for finding a logical measure of progress for planning. (McCarthy 1982) reiterates a reduction alluded to by (Kempe 1879 [1]) that will never require backtracking: if a country $C$ in a map has three or fewer neighbors, then we can postpone four-coloring $C$ until the other three neighbors have been colored. The original problem of coloring reduces to the same map minus $C$. In some cases, such as the map of the United States, this process will continue until the map is completely stripped, in which case coloring is done in the reverse order that the countries are eliminated, respecting other countries' colors as required. She provide a theorem which elucidates the fundamental structure required for maps to be so easily reduced, and shows how they are in fact, antimatroids. The shells removed are the vertices of degree three or less and generalize to $n$ colors:

***Definition 1*** ($n$-reducible maps): *A map* $(V,E)$ *is* $n$-reducible *if one can repeatedly remove vertices of degree* $n$ *or less from the graph, until the empty graph is encountered.*

If a map $(V,E)$ is $n$-reducible then we can color it with $n+1$ colors without ever having to backtrack.

***Theorem 1*** ($n$-reducible maps are antimatroids)[14]**:** *Let* $(V,E)$ *represent a map*. $(V,E)$ *is* $n$-*reducible iff* $(V,\mathcal{L}H)$ *is an antimatroid*.

There are many equivalent formulations of the four color problem e.g., [17],[18] and in fact in the absence of strong hopeful attempts for some notorious conjectures a quite strange method have been developed for possible solving time[3]. Unfortunately Parmar's idea for *n*-reducible maps is not good enough to show that all maps are four-colorable (without use of an computer). Here we give a new idea which is an antimatroid and pre-planning of the coloring algorithm is enough to show that all maps are four colorable. As in [1],[2] we assume that *G* is an internally 5-connected maximal planar graph.

[3] In SIGACT News Complexity Theory Column 36, Guest Column: The P=?NP Poll, William I. Gasarch have taken a poll of theorists to see what they think about P=?NP and even one person said that the conjecture will be settled like the proof of the four color theorem!

# 3 The need of an innovative approach and spiral chains

Throughout the paper by *G* we assume an *n* vertex maximal planar graph drawn on the plane so that no two edges cross. We only note that in *G* all faces (smallest cycles are triangles) including exterior cycles are triangles. We obtain spiral chain(s) of *G* as follow (see Fig.1). Let $v_a, v_b, v_c$ be the three vertices of the outer triangle of *G*. We obtain the spiral chain by starting from $v_a$ and selecting the edges of the outer triangle in the clockwise direction (through the vertex $v_c$) and enter in the inner vertex of *G* from $v_b$ by selecting the left-most edge (not the outer triangle edge $\{v_b, v_a\}$) going into *G*. The subgraph $G_{i-1}$ of *G* is obtained by deleting the vertices $v_a, v_b, v_c$ of the outer triangle. Then repeat the same operation for $G_{i-1}$ and others till scan all vertices of G. If there is no edge connecting the last vertex of the current partial spiral chain to the vertex of $G_j$ for some *j* then call the chain the first spiral chain $S_1$ and repeat the whole operation described above (starting from the closest vertex to the last vertex of $S_1$ in the clockwise direction) for the vertices of $G_j$ . Eventually we obtain vertex disjoint spiral chains $S_1, S_2, ..., S_k, k \geq 1$ so that $\cup_{i=1}^{k} |V(S_i)| = n$ . Note that if $k = 1$ then $S_1$ is also an Hamiltonian path of *G*.

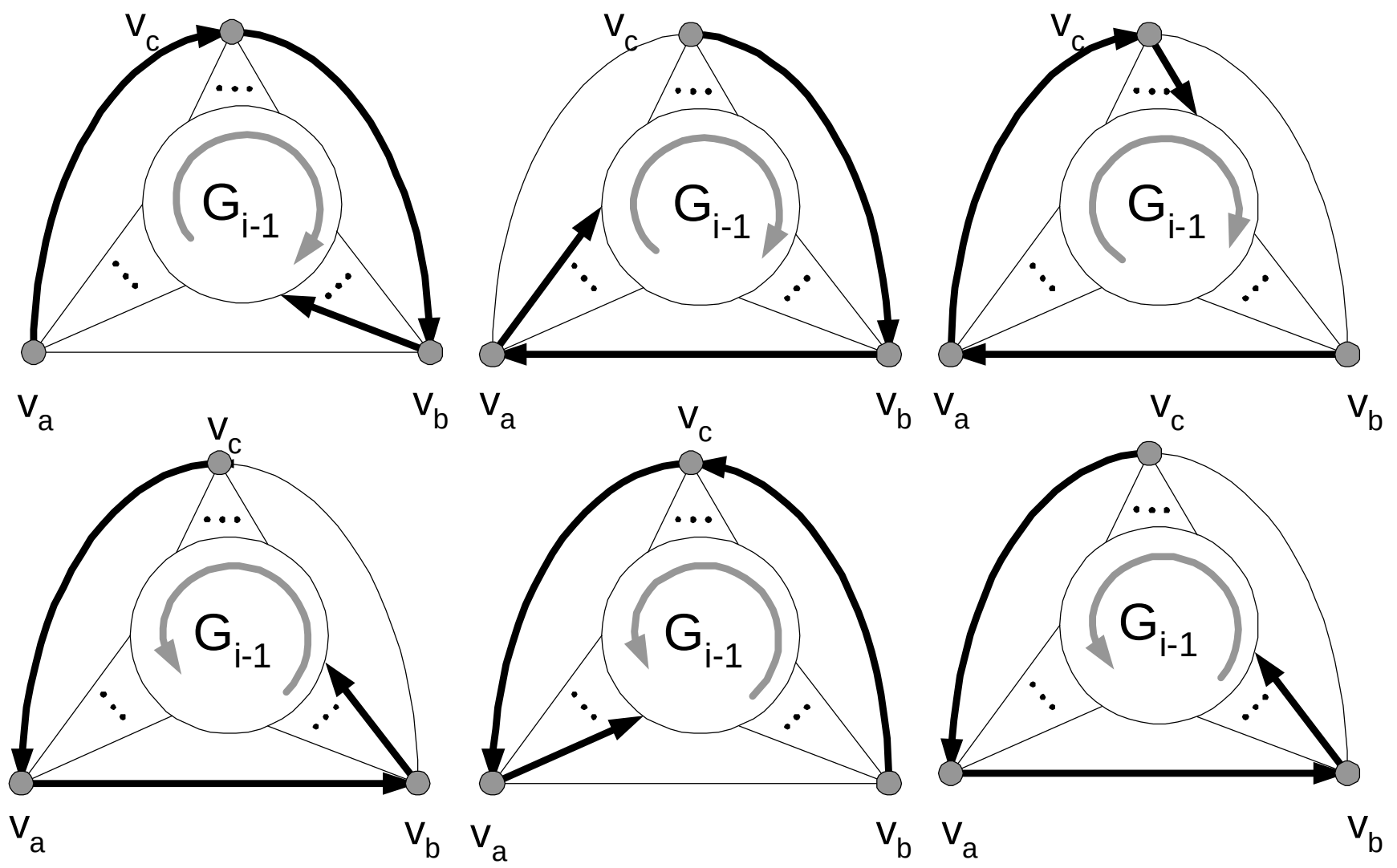

Fig. 1

In general a vertex in the spiral chain looks like a way to "cut up" a graph so that each vertex is connected to one of four kinds of vertices: one vertex forward of it in the spiral, and one vertex behind it, and then a set of vertices to its "right" which are bisected by the chain, and another set on the "left" which are bisected on the other side. Consider now any pairwise neighbor three sub-spiral chains $S_{i+1}, S_i, S_{i-1}$ in *G* such that vertices of $S_{i-1}$ (right sub-chain) have already been colored completely (say, with 3,2,4,3,2 in Fig. 2(a)), vertices of $S_i$ (middle sub-chain) colored partially (say 1,2,4 in Fig. 2(a)) till the vertex *x* and all vertices of $S_{i-1}$ (left sub-chain) have not yet been colored. In the spiral chain coloring algorithm in deciding the color $c \in \{1,2,3,4\}$ of the vertex $x$ we need only to consider the color of the vertex $y$ on $S_i$ with $(y,x) \in S_i$ and the color of the vertices adjacent at $x$ on $S_{i-1}$. In another words while coloring the current vertex $x$ we only to consider the vertices of right-spiral chain and the color of the vertex $y$.

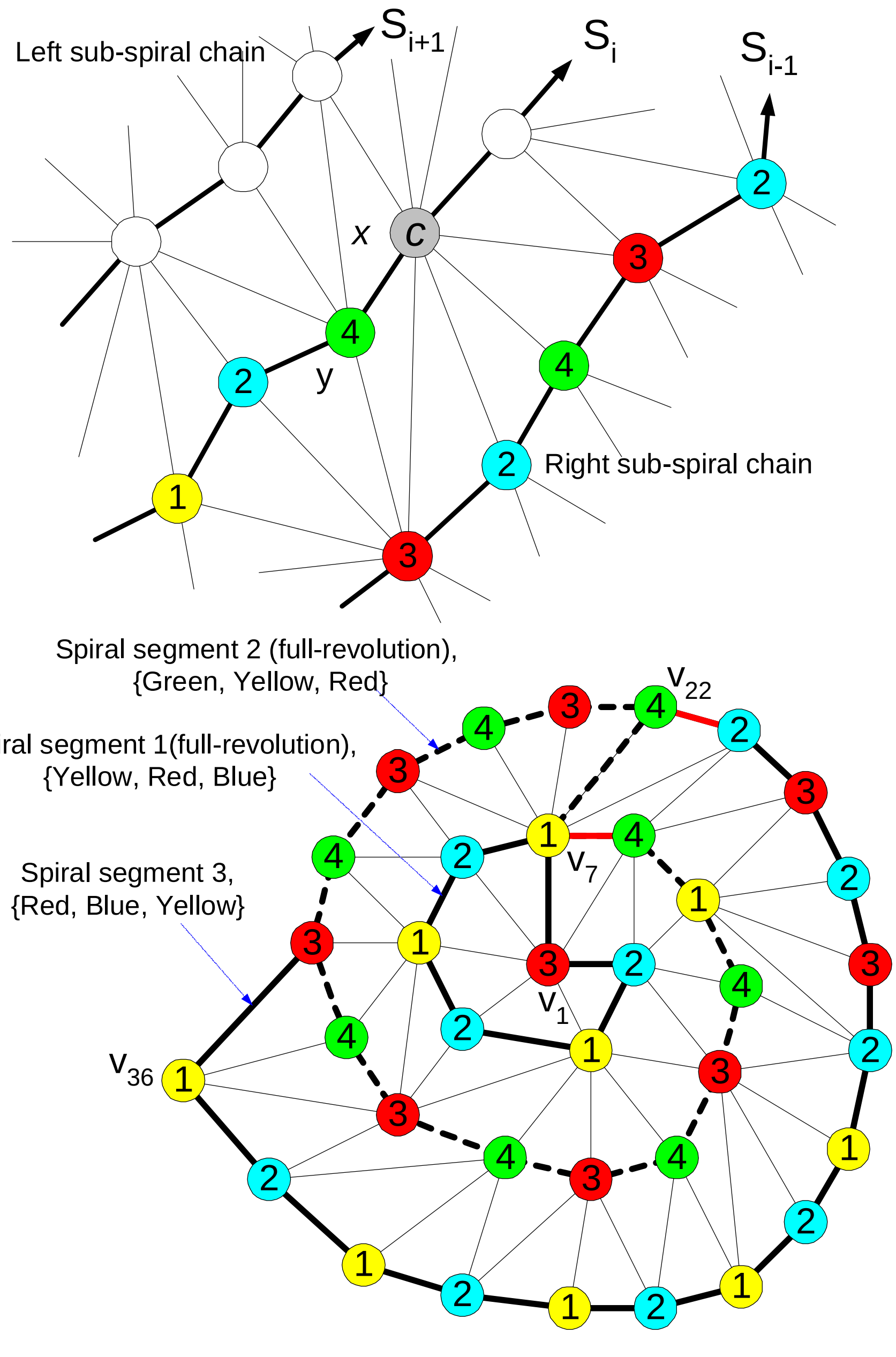

Fig. 2(a) and (b)

Another useful concept to be used later in the spiral chain coloring, is the segmentation of a spiral chain.

**Definition 2.** *Let* $S_n$ *be a n-vertex spiral chain with the ordered vertex set* $V=\{v_1,v_2,v_3,...,v_n\}$. *We say the subset of consecutive vertices* $V\langle i,j\rangle=\{v_i,v_{i+1},v_{i+2},...,v_j\}$, *with* $i\geq 1$ *and* $i<j\leq n$ *of vertex set* $V$ *is a spiral-segment with full revolution if* $v_j$ *is the last vertex in* $S_n$ *such that* $(v_j,v_i)\in E(G)$ *and* $(v_j,v_{i+1})\notin E(G)$.

As an example to spiral-segments of the triangulated planar graph with *36* vertices shown in Fig. 2(b) are $V\langle 1,7\rangle, V\langle 8,22\rangle$ and $V\langle 23,36\rangle$. It is clear that if the graph $G$ has only one spiral-chain then it can be partioned into spiral segments as

$$V=\{V\langle 1,i_1\rangle\cup V\langle i_1+1,i_2\rangle\cup...\cup V\langle i_{(k-1)}+1,i_k\rangle\}.$$

for some $k\geq 1, i_k=n$. If $G$ has more than one spiral-chain, say $p\geq 2$ chains, then we can generalize the above as

$$V=\{\cup_{i=1}^{k_1}V_1\langle i,i+1\rangle\cup_{i=k_1}^{k_2}V_2\langle i,i+1\rangle\cup...\cup_{i=k_{p-1}}^{k_p}V_p\langle i,i+1\rangle\}.$$

For any two consecutive spiral segments $V\langle i_1,i_2\rangle$ and $V\langle i_2+1,i_3\rangle$ we say that the second one surrounds the first one. Let $C_1$ and $C_2$ be the 3-coloring of spiral segments $V\langle i_1,i_2\rangle$ and $V\langle i_2,i_3\rangle$. If the color classes of $C_1$ and $C_2$ differ at least by two then we say that the 3-coloring of the consecutive spiral segments are *safe* (has a safe 3-coloring). In fact safety of consecutive spiral-segments is crucial in our proof of the main theorem since it gives us relaxation for the proper coloring of the whole maximal planar graph *G*.

## 4 Some properties of spiral chains

A *theta* graph is a homeomorph of $K_{2,3}$ or equivalently, a pair of distinct vertices joined by three pairwise internally-vertex disjoint paths. A generalized theta graph is a spiral-chain separating graph which is an maximal outerplanar subgraph $G_\theta$ of $G$ between any two consecutive spiral chains $S_i$ and $S_{i+1}$ such that vertices $x\in S_i, y\in S_{i+1}$ and the edge $(x,y)\notin G_\theta$. Note that there is only one path from *x* to *y* through the vertices of $S_i$, hence $(x,y)\notin G_\theta$ preventing to join $S_i$ and $S_{i+1}$ into a single spiral chain. In the Fig. 3(a) we have illustrated this that spiral chain $S_1$ (blue chain) cannot continue its way, i.e., to vertex y of the spiral chain $S_2$ (red chain) from vertex *x* to *y* because of the theta separator subgraph $G_\theta$. In Fig. 3(b) we have another situation in which in spiral chain $S_1$ the edge *(x,y)* must be chosen instead of the edge *(x,t)* connecting spiral chain $S_2$ which results two spiral chains. Note also that if the size of spiral chain $|S_2|=1$ then $S_2$ reduces to a single vertex *y*.

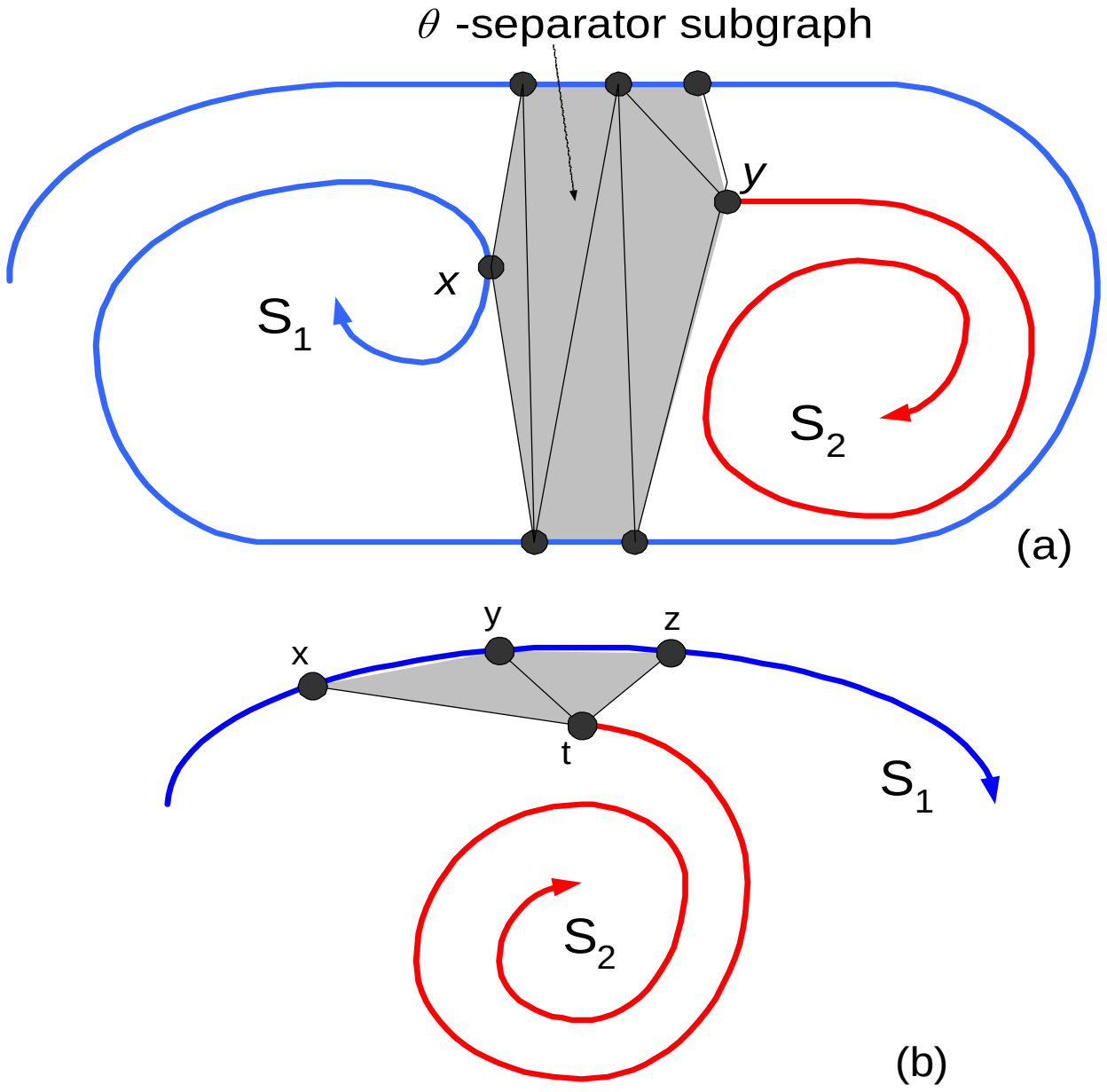

Fig. 3(a) and (b)

A *fan* $F_{n-2}$ (here $(n-2)$ is the number of triangles) is a graph obtained from a *n*-vertex path $P_n$ by joining all its vertices to another vertex (apex) $v_{n+1}$ or equivalently can be obtained by deleting any outer edge $(v_i, v_{i+1})$ of an $(n+1)$-vertex wheel $W_{n+1}$. Ladder-fan graph $G_l$ is obtained by the union of *k* fans such that all its vertices lie either on one of the two parallel lines $L_u$ (upper-line) or $L_l$ (lower-line). We call a fan $F(u)_i \in L(f)$ (respectively, $F(l)_i \in L(f)$) an upper-fan (lower-fan) of $G_l$ if its apex-vertex is on $L_u$ (respectively, if its apex-vertex is on $L_l$). For example $F(u)_2 \cup F(l)_3 \cup F(u)_4 \cup F(l)_2$ is ladder-fan $G_l$ consisting of, respectively with *2* upper, *3* lower, *4* upper and *2* lower triangles. In a fan decomposition of the spiral chains if all vertices of the fan lie on the same side upper or lower-line, then we call such a fan *c*-type (or "*corner*" type) fan. As an example in Fig.6 triangles numbered as *1*, *3* and *4* are *c*-type fans (triangles). Actually this situation arise when the spiral chain *S* cuts a wheel into two through two spoke-edges. In this case we have to switch from one 3-color class into another e.g., from {blue, green, yellow} to {red, green, yellow} for the fans *3* and *4* in Fig.7.

***Lemma 1:*** *A ladder-fan graph* $G_l$ *is 3-colorable.*

*Proof:* It follows from the fact that $G_l$ is a maximal outerplanar and any triangle in $G_l$, with the exception that of the left-most and right-most triangles, share an edge between the two parallel lines (paths) $L_u$ and $L_l$. □

***Lemma 2***: *Two spiral chains* $S_i$ *and* $S_j$, $i \neq j$ *are separated only by a theta separator subgraph* $G_\theta$.

*Proof:* From Fig. 3(a) and (b) we see that the edge(s) of $S_1$(theta subgraphs shown in gray) blocks the spiral chain $S_1$ to join with the spiral chain $S_2$. Otherwise we will have only one spiral chain e.g., delete the edge $(x, y)$ in Fig. 3(b). □

***Lemma 3***: *Let $S_i$ and $S_{i+1}$ be the two consecutive spiral segments and let $G_{i,i+1}$ be the ladder-fan formed by these segments. Then there exists four coloring of $G_{i,i+1}$ consisting union of 3-colorings of the spiral-segments so that its color classes differ by two colors.*

*Proof:* Since $G_{i,i+1}$ is an outerplanar subgraph of $G$ there exists at least 3-coloring of the ladder-fan. Let $C_i = \{1,2,3\}$ and $C_{i+1} = \{2,3,4\}$ be the color classes of the consecutive spiral segments $S_i$ and $S_{i+1}$. That is the label (color) *1* is not in $C_{i+1}$ and *4* is not in $C_i$. Assume first that in $G_{i,i+1}$ there is no edge $e = \{(v_p, v_q) \mid v_p, v_q \notin S_{i+1}, |p-q| \geq 2\}$ between any two non-consecutive vertices of $S_i$ or $S_{i+1}$. We consider two cases:

(a) $v_a \in S_i$ with color $f(v_a) = \{2\}$ or $\{3\}$ are adjacent to $v_j \in S_{i+1}, j = 1, 2, ..., k$ vertices.

(b) $v_a \in S_i$ with color $f(v_a) = \{1\}$ are adjacent to $v_j \in S_{i+1}, j = 1, 2, ..., k$ vertices.

In the case (a) if $f(v_a) = 2$ then color alternatingly $k$ adjacent vertices on $S_{i+1}$ with *3* and *4* and if $f(v_a) = 3$ then color $k$ adjacent vertices with colors *2* and *4*. The case (b) is even simpler than the case (a) since $C_{i+1} \cap \{1\} = \varnothing$. Therefore in worst we will have 3-colorings for each spiral segments $S_i$ and $S_{i+1}$. Next assume that we have non-consecutive edges $(v_i, v_{j+2}), (v_{i+2}, v_{i+4}), (v_i, v_{i+4}) \in E(G)$ where $v_i, v_{i+2}, v_{i+4} \in V(S_{i+1})$. We have shown 3-coloring of $S_{i+1}$ under the possible color assignment of the vertex $v_a \in S_i$ in Fig. 5. Note that if we would insert the edge $(v_{i+1}, v_{i+3}) \in E(G)$ then by definition of the spiral chain the topology of $S_i$ and $S_{i+1}$ would be changed to another spiral chains $S_i'$ and $S_{i+1}'$ and we then re-color $S_{i+1}'$ with the colors of $C_{i+1} = \{2,3,4\}$ again. □

Another way of proving the above lemma is to note that the induced planar graph $G_{i,i+1}$ on $V(S_i \cup S_{i+1})$ is an maximal outerplanar. Now first color the vertices of $G_{i,i+1}$ with *1,2,3* and then change all vertices with color *1* in $S_{i+1}$ to color *4*. Now we have $C_i = \{1,2,3\}$ and $C_{i+1} = \{2,3,4\}$.

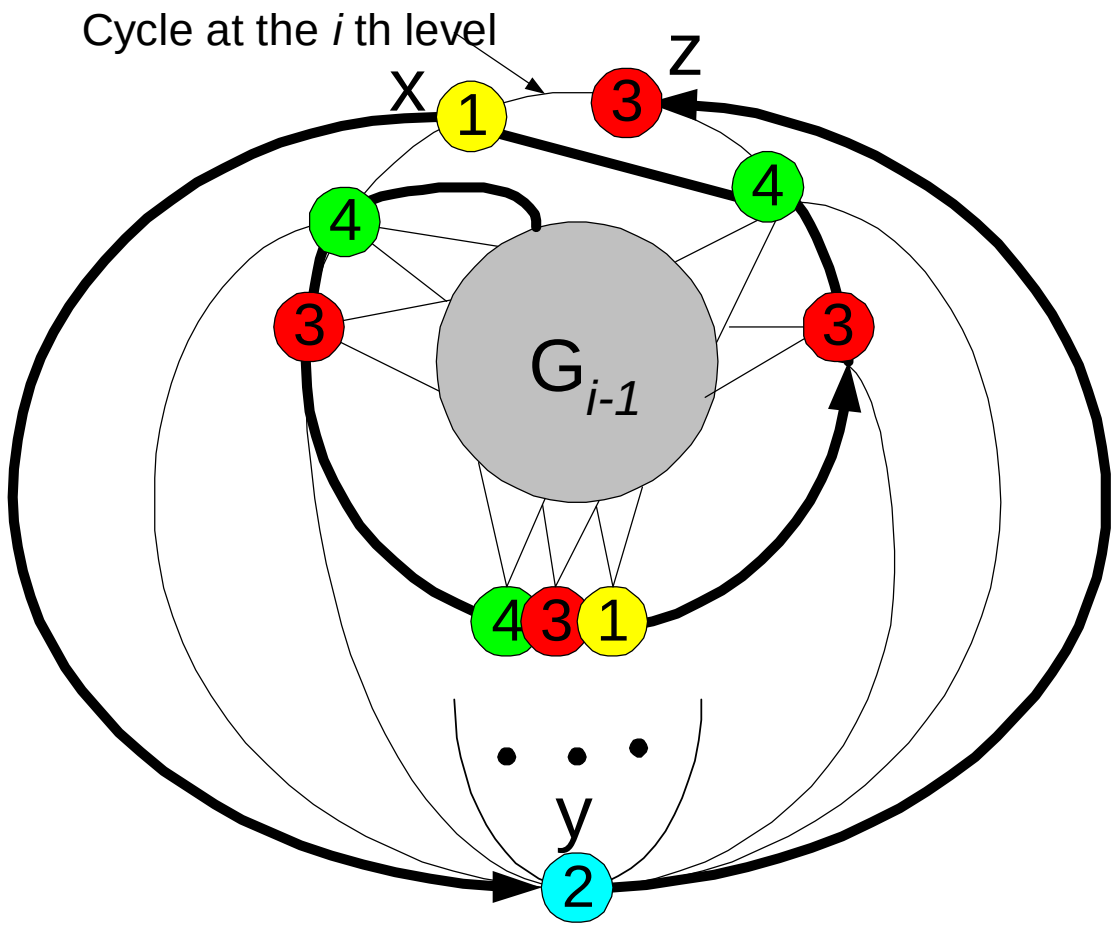

Coloring outercycle with the spiral chain
(termination condition)

Fig. 4

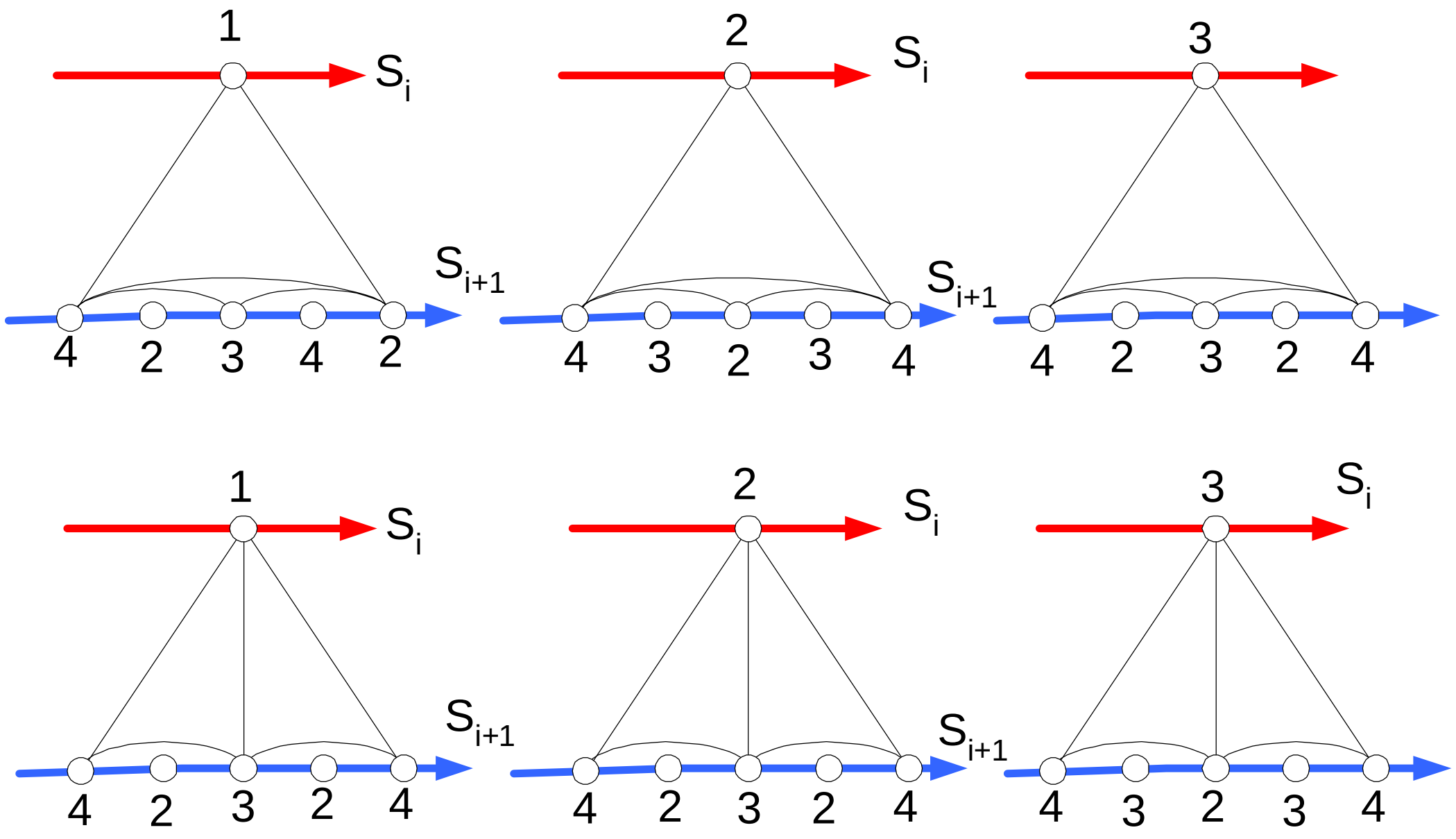

Fig. 5

## 5 Four color of the historical three "bad" examples

In Fig. 6 we have shown the famous counter-example graphs to the Kempe's method. The figures also illustrate the use of spiral chains in the four coloring of these graphs. Fig. 7 is another maximal planar graph from [8] for an illustration spiral chain four-coloring as well as its fan-decomposition.

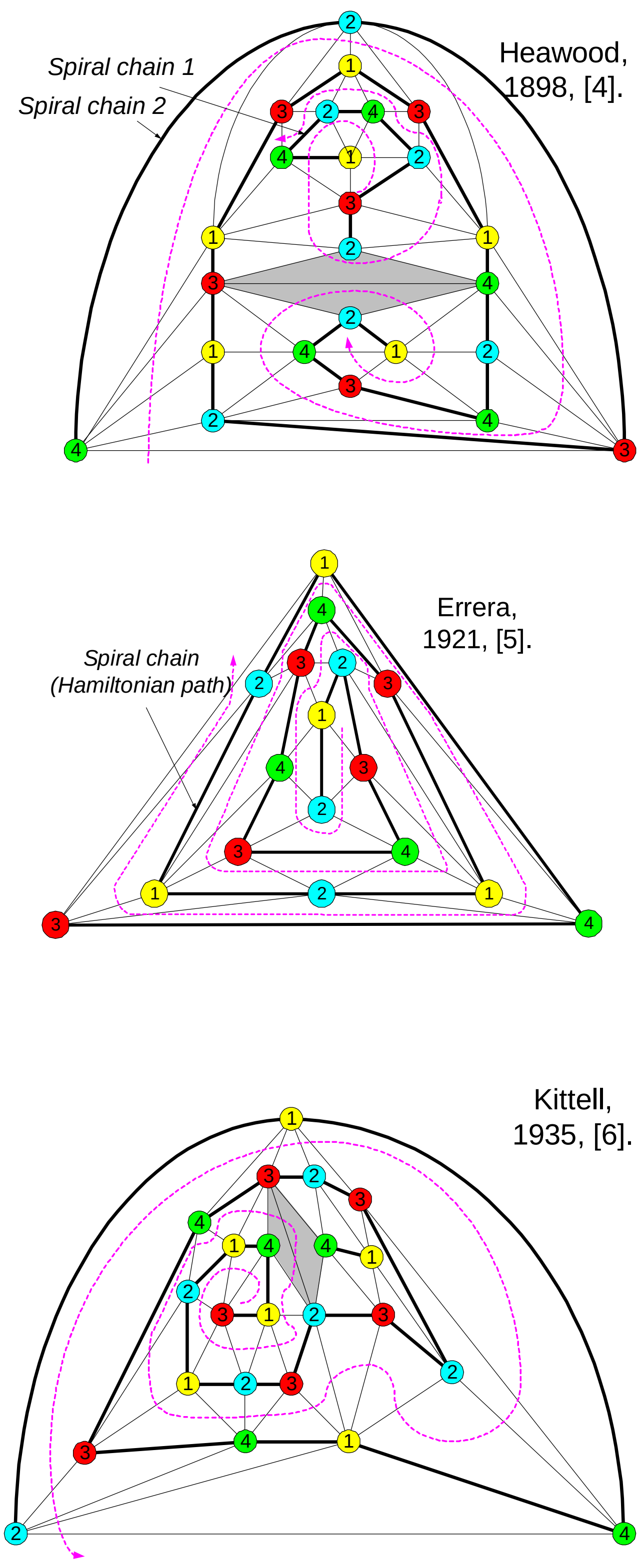

Fig. 6

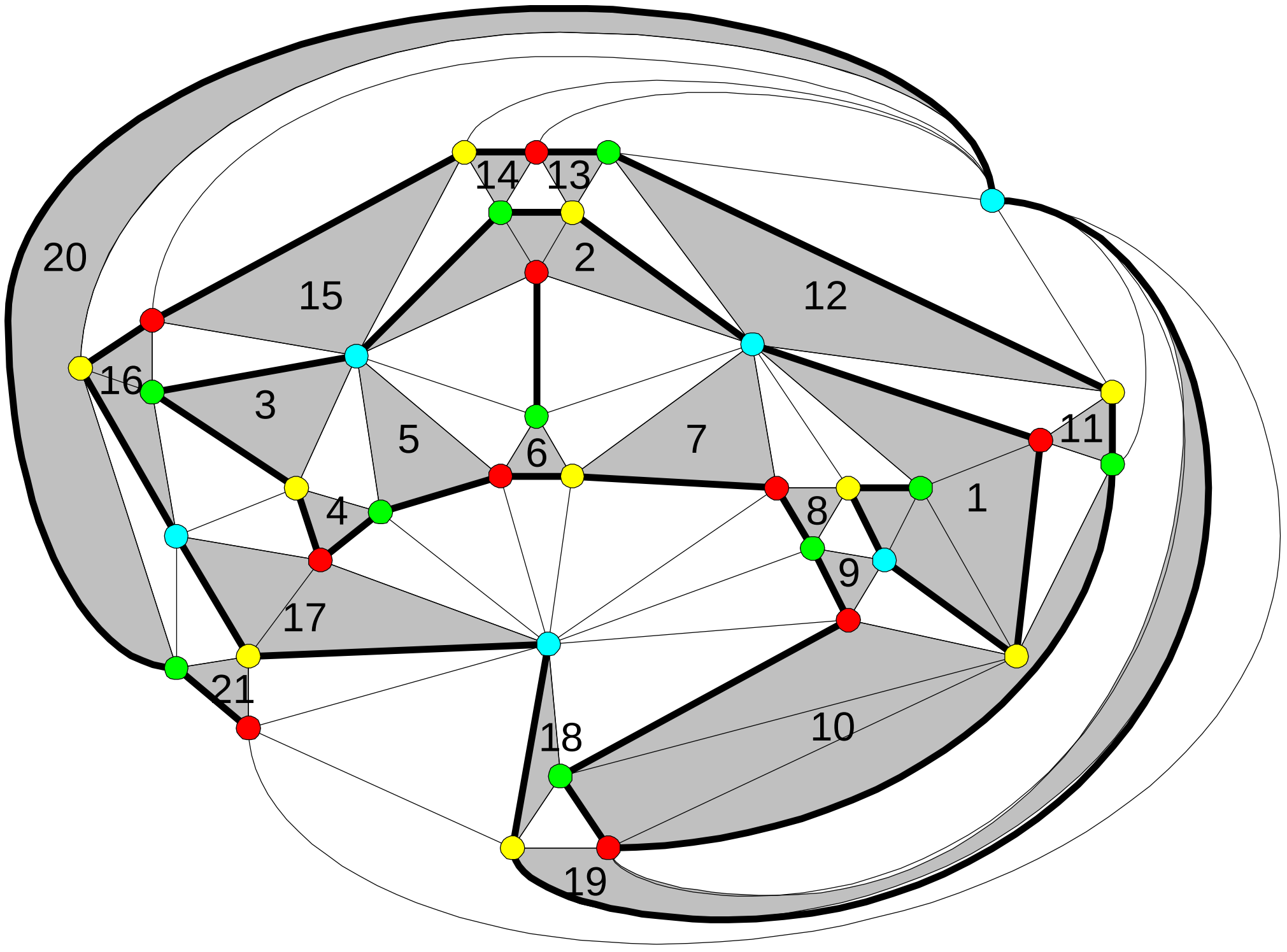

Fig. 7

# 6 Spiral-chain coloring algorithm

***Lemma 4***: *Spiral chains of a maximal planar graph enables a fan-decomposition of G.*

*Proof:* Details of this lemma is left to the reader (see Fig. 7).

***Theorem 2*:** Spiral chains coloring are shelling structures.
*Proof:* Let $S_1, S_2, ..., S_k$ be the spiral chains of a maximal planar graph *G*. Let $S_i, 1 \leq i \leq k$ be further partioned into $k_i$ spiral-segments $Si = \{S_{i_1}, S_{i_2}, ..., S_{i_{k_i}}\}$. From Lemma 3 we know that for any consecutive spiral segment there exists 3-coloring for which color classes differ by at least *2*. Now starting from the inner spiral chain (i.e., reverse order of the spiral chains $S_1, S_2, ..., S_k$) and color each spiral-segment with the suitable three colors (Lemma 3). Repeating this process one by one for all spiral chains clearly results in a four coloring of the graph *G*. Perhaps we would say a little bit more for the last step (termination condition). Coloring the last spiral-segment of $S_1$ must be properly colored with the three colors. In Fig. 4 the last spiral-segment has been colored by 1, 2, and 4 (yellow, blue, red) while the neighbor spiral segment colored by 1, 3, 4 (yellow, red, green). Therefore four coloring of G is completed. Now by deleting all colored spiral-segments we get an empty graph. □

Detailed discussions and pervious results derived leads us to state the following main theorem of the paper:

***Theorem 3***: *Maximal planar graphs are four colorable.*

***Corollary 1***: *All planar maps are four colorable.*

# 7 Conclusion

We have given a new solution of the four color theorem by the use of spiral chains in the maximal planar graphs. We concept of spiral chain is new and quite powerful and appears to be useful for some other problems in topological graph theory.

August 2004